\documentclass[12pt]{amsart}
\usepackage{amssymb,latexsym}
\usepackage[all]{xy}
\newtheorem{thm}{Theorem}[section]
\newtheorem*{thm*}{Theorem}
\newtheorem{proposition}[thm]{Proposition}
\newtheorem{lemma}[thm]{Lemma}
\newtheorem{cor}[thm]{Corollary}
\theoremstyle{definition}

\newtheorem{notation}[thm]{Notation}
\newtheorem{remark}[thm]{Remark}

\renewcommand{\deg}{{\rm deg}\,}
\newcommand{\by}[1]{\stackrel{#1}{\longrightarrow}}
\newcommand{\Pic}{{\rm Pic}\,}

\newcommand{\End}{{\rm End}\,}

\newcommand{\sEnd}{\mbox{${\sE}{nd}$}}

\newcommand{\boxtensor}{{\Box\kern-9.03pt\raise1.42pt\hbox{$\times$}}}
\newcommand{\C}{{\mathbb C}}
\newcommand{\Z}{{\mathbb Z}}
\newcommand{\Hm}{{\mathbb H}}
\newcommand{\V}{{\mathbb V}}
\newcommand{\U}{{\mathbb U}}

\newcommand{\W}{{\mathbb W}}
\newcommand{\tensor}{\otimes}

\newcommand{\bP}{{\mathbf P}}
\newcommand{\sE}{{\mathcal E}}
\newcommand{\sI}{{\mathcal I}}
\newcommand{\sF}{{\mathcal F}}

\newcommand{\sL}{{\mathcal L}}
\newcommand{\sM}{{\mathcal M}}
\newcommand{\sN}{{\mathcal N}}
\newcommand{\sO}{{\mathcal O}}

\newcommand{\sT}{{\mathcal T}}

\renewcommand{\tilde}{\widetilde}

\numberwithin{equation}{section}
\newcounter{elno}                

\newcounter{example}[section]
\def\theexample{\thesection.\arabic{example}}

\begin{document}

\title{Restriction theorems for homogeneous bundles}
{}

\author{V.B. Mehta and V. Trivedi}
\address{School of Mathematics, Tata Institute of
Fundamental Research,
Homi Bhabha Road, Mumbai-400005, India}
\email{vikram@math.tifr.res.in and vija@math.tifr.res.in}

\subjclass{14L30}
\date{}
\begin{abstract}{We prove that for an irreducible representation
$\tau:GL(n)\to GL(W)$, the associated
homogeneous $\bP_k^n$-vector
bundle $\W_{\tau}$
is strongly semistable when restricted to any smooth quadric or to any
smooth cubic in $\bP_k^n$, where $k$ is an algebraically closed field of
characteristic $\neq 2,3$ respectively. 
In  
particular $\W_{\tau}$ is  semistable when restricted to 
general 
hypersurfaces of degree $\geq 2$  and is strongly semistable when
restricted
to the 
$k$-generic hypersurface of  degree $\geq 2$.}
\end{abstract}
\maketitle
\section{Introduction}
In this paper we study the semistable restriction theorem for the
homogeneous vector
bundles on $\bP^n_k$ which come  from irreducible $GL(n)$-representations.

In general suppose $G$ is  a reductive algebraic
group over an algebraically closed field $k$ and  $P\subset G$ is  a
parabolic group. Then there is an equivalence between the category of 
 homogeneous $G$-bundles over $G/P$ and 
the category  of $P$-representations, where
a $P$-representation $\rho: P\to GL(V)$ on a
$k$-vector space 
$V$ induces a homogeneous $G$-bundle
$\V_{\rho}$ on $G/P$ given by 

$$\V_{\rho} = \frac{G\times V}{P} = 
\frac{G\times V}{\{(g,v)\cong
(gh,h^{-1}v)\mid
g\in G, v\in V, h\in P\}}.$$

Now for the rest of the paper we fix the following 
\begin{notation}\label{n1} The field $k$ is an algebraically closed
field and $G = SL(n+1,k)$, and  $P$
is 
the maximal parabolic subgroup of $G$ given by
$$ P = \left\{\left[\begin{array}{cc}
g_{11} & *\\
0 & A \end{array}\right]
\in SL(n+1), ~~\mbox{where}~~A\in GL(n)\right\} $$
and $G/P\simeq \bP_k^n$ is a  canonical isomorphism.
\end{notation}

Now, if $\sigma:GL(n)\to GL(V)$ is an irreducible
$GL(n)$-representation then it induces an irreducible $P$-representation
$\rho:P\to GL(V)$ given by 

\begin{equation}\label{e1}\left[\begin{array}{cc}
g_{11} & *\\
0 & A \end{array}\right] \mapsto \sigma(A),
\end{equation}

\noindent which gives a $G$-homogeneous bundle on $G/P = {\mathbf P}_k^n$.
Conversely, any $G$-homogeneous bundle $\V$, given by an irreducible
$P$-representation
$\rho:P\to GL(V)$, is in fact induced by an irreducible
$GL(n)$-representation (upto tensoring by $\sO_{\bP^n_k}(r)$,
 for some $r$). 

In this paper we prove the following

\begin{thm}\label{1}Let  
$\tau:GL(n)\to GL(W)$ be an irreducible
$GL(n)$-representation, where $W$ is a $k$-vector space. Let
$\W_{\tau}$ be the
associated $G$-homogeneous bundle on $G/P = \bP_k^n$. Let 
\begin{enumerate} 
\item $X$ = smooth quadric, if ${\rm char}~k \neq 2$, or
\item $X$ = smooth  cubic, if ${\rm char}~k \neq 3$,
\end{enumerate}
Then the
bundle $\W_{\tau}\mid_X$ is
strongly semistable.\end{thm}

We note that Theorem~\ref{1}  implies  $\W_{\tau}$ itself is
 semistable on $\bP^n_k$.
 However this result, in much more general form, has been 
proved in  [R], [U], [MR1] and [B].

Theorem~\ref{1} implies (see Corollory~\ref{c1}) that, provided ${\rm
char}~k \neq 2, 3$, 
 the bundle $\W_{\tau}\mid_H$ is
{\it semistable}, for a general
hypersurface $H$ 
 of  degree $\geq 2$ in $\bP_k^n$,  
and $\W_{\tau}\mid_{H_0}$ is 
{\it
strongly
semistable} for generic hypersurface $H_0$ of degree $d\geq 2$.
This is equivalent to the statement that, given  $s\geq 0$, the $s^{th}$
Frobenius pull back $F^{s*}\W_{\tau}\mid_H$ is semistable for a general
hypersurface $H$ of  degree $\geq 2$ in $\bP^n_k$. 
Moreover when the bundle $\W_{\tau}$ comes from the standard
representation, {\it i.e.}, $\W_{\tau}$ is the tangent bundle (upto a
twist
by
a line bundle) of $\bP^n_k$, where $n\geq 4$, then we can prove a stronger
statement, by
replacing
the word `semistable' by `stable' everywhere in Theorem~\ref{1} and
Corollory~\ref{c1}.

In this context we recall that, Mehta-Ramanathan~[MR2]
have proved that
if $E$ is a semistable sheaf on a smooth projective variety (over
a field of arbitrary characteristic) then $E$
restricted to a general hypersurface  of degree $a$ (where $a$ is
any sufficiently large integer) is semistable. On the other hand, Flenner
[F] proved this assertion, where the degree $a$ of the hypersurface
depends only on the rank of $E$ and degree of the variety $X$, provided
the characteristic is $0$.

The paper is organised as follows: In Section~2, we recall some general
facts about smooth quadrics. Then we discuss the vector bundle
$\V_{\sigma} = \sT_{\bP^n_k}(-1)$ associated to the standard
representation
$\sigma:GL(n)\by{} GL(V)$
 and its restriction to smooth quadrics. In particular, for a smooth
quadric $Q \subset \bP^n_k$, we show
that
$\V_{\sigma}\mid_Q$ has a unique $SO(n+1)$-homogeneous proper subbundle,
if
$n\geq 4$, 
(see remark~\ref{r2} for details).

In
Section~3, we prove that if char~$k \neq 2$ then  $\sT_{\bP^n_k}\mid_Q$ 
is strongly stable  if $n\geq 3$, and is strongly semistable if $n =2$. 
Moreover the tangent bundle $\sT_Q$ of $Q$ is semistable and is of  
positive
slope. 

In Section~4 we prove that, if char~$k\neq 3$ and 
 $X \subset \bP^n_k$ is an arbitrary smooth cubic
hypersurface then
$\sT_{\bP^n_k}\mid_X$
is strongly stable if $n\geq 4$ and strongly semistable if $n=2$ or $n=3$. 
 Moreover the tangent bundle $T_X$ of $X$ is either stable if
$n\neq 3$, or $\mu_{\rm min}(\sT_X) \geq 0$ if $n=3$.
In fact, we show that the argument given in [PW],
 to prove stablity of $\sT_X$, for a smooth hypersurface
of $\deg d \geq 3$, $n\geq
4$ and $k =\C$, can be modified so as to work
over any algebraically closed field of characteristic coprime to  $d$
(this hypothesis is needed so that the cup product 
with $c_1(\sO_{\bP^n_k}(d))$ is an injective map). 

Finally in Section~5, we show (see Theorem~\ref{1}) that, if
$\V_{\sigma}\mid_X$ is semistable and
$\mu_{\rm{min}}(\V_{\sigma}\mid_X) \geq
0$, where $X$ is a smooth hypersurface in $\bP^n_k$
then the bundle $\W_{\tau}\mid_X$ is strongly semistable for any
irreducible representation $\tau:GL(n)\by{} GL(W)$.

 \section{Some general facts about quadrics}
 \subsection{Embedding of
quadrics in
${\bP}_k^n$}
 Let $V$ be a vector-space of
dimension $n+1$ over $k$
(characteristic $k \neq 2$).
Let us choose a basis $\{e_1,\ldots, e_{n+1}\}$ of
$V$. Represent a point $v\in V$ by 
$$ v = (x_1,\ldots, x_{n/2}, z, y_1, \ldots,
y_{n/2}),~~\mbox{if}~~ n ~~\mbox{is even},$$ 
$$ v = (x_1,\ldots, x_{(n+1)/2}, y_1, \ldots,
y_{(n+1)/2}),~~\mbox{if}~~ n~~\mbox{ is odd},$$ 
with respect to the basis $\{e_1,\ldots, e_{n+1}\}$.
Without loss of generality, one can assume that any
fixed smooth quadric $ Q\subset \bP_k^n $ is given by the quadratic
form
$$ {\tilde Q(v)} = z^2 + 2(x_1y_{n/2} + \cdots +
x_{n/2}y_1),~~\mbox{~~if}~~
n~~~ \mbox{is even and} $$
$${\tilde Q(v)} = x_1y_{(n+1)/2} + \cdots +
x_{(n+1)/2}y_1),~~\mbox{~~if}~~
n ~~\mbox{is odd}.$$

Let $$ \begin{array}{lcl}
SO(n+1) & = & \{A\in SL(n+1)\mid {\tilde Q(Av)} = {\tilde Q(v)}~~\mbox{for
all}~~v\in V\}\\
 & = & \{A\in SL(n+1) \mid A^tJA = J\}\end{array},$$
where 
$$ J = \left[\begin{array}{ccc}
0  & \cdots & 1\\
0  & { \begin{array}{ccc}
{}  & {} & \cdot\\
{} &  \cdot & {}\\
\cdot & {} & {} \end{array}} & 0\\
1  & \cdots & 0 \end{array}\right] \in GL(n+1).$$

\begin{notation}\label{n2}Let  $P_1 = P \cap SO(n+1)$ denote the maximal
parabolic group in $ SO(n+1)$
such that 
$$\displaystyle{\left\{\left[\begin{array}{ccc}
a_{11} & 0 &  0 \\
0 & A &  0 \\
0 & 0 & a_{11}^{-1}\end{array}\right],~~\mbox{where}~~ A\in
SO(n-1), a_{11}\in
k^*\right\}} \subseteq P_1,~~\mbox{and}~~$$
$$P_1 \subseteq 
\displaystyle{\left\{\left[\begin{array}{ccc}
a_{11} & * & *\\
0 & A & *\\
0 & 0 & a_{11}^{-1}\end{array}\right],~~\mbox{where}~~ A\in SO(n-1),
a_{11}\in
k^*\right\}}.$$\end{notation}
Then we have the canonical identification
$$\begin{array}{ccc}
\bP_k^n &  \simeq & SL(n+1)/P\\
\uparrow & & \uparrow\\
Q & \simeq & SO(n+1)/P_1.
\end{array}$$

\vspace{10pt}

\subsection{Standard representation of $GL(n)$}

Consider the canonical short exact sequence of sheaves of
$\sO_{{\bP}_k^n}$-modules
$$0 \by{} \Omega^1_{\bP_k^n}(1) \by{} H^0(\bP_k^n,
\sO_{\bP_k^n})\tensor{\sO_{\bP_k^n}} \by{}
\sO_{\bP_k^n}(1)\by{} 0.$$
The dual sequence is
\begin{equation}\label{e2}
0 \by{} \sO_{\bP_k^n}(-1) \by{} H^0(\bP_k^n,
\sO_{\bP_k^n})\tensor{\sO_{\bP_k^n}} 
\by{}
\sT_{\bP_k^n}(-1) \by{} 0,\end{equation}
where $\sT_{\bP^n_k}$ is the tangent sheaf of $\bP_k^n$. Now this sequence
is
also a short exact sequence of $G$-homogeneous bundles on $G/P =
\bP_k^n$ (see \ref{n1}). Hence
there exists a corresponding short exact sequence of $P$-modules 
$$0\by{} V_{2} \by{f} V_{1} \by{\eta} V\by{} 0,$$
where the $P$-module  structure is given as follows.

 Let $V_1$, $V$ and
$V_2$
be $n+1$, $n$ and $1$ dimensional $k$-vector spaces respectively,
with fixed bases. Let 
$f: (c) \mapsto (c, 0,\ldots,0)$ and let
$$\eta:(a_1,\ldots, a_{n+1})\mapsto (0,a_2,\ldots, a_{n+1}).$$

Now representing the
elements
of the vector spaces as coloumn vectors and expressing any $g\in P$ as
$$ g = \left[ \begin{array}{cc}
g_{11} & *\\
0 & B\end{array} \right], ~~~\mbox{where}~~B\in GL(n),$$
we define the representations as follows:

The representation $\rho_1:P\by{} GL(V_1)$ is given by  
$$\rho_1(g)\left[\begin{array}{c}
a_1\\
\vdots\\
a_{n+1}\end{array}\right] = [g]\left[\begin{array}{c}
a_1\\
\vdots\\
a_{n+1}\end{array}\right].$$
The representation $\rho_2: P\by{} GL(V_2)$ is given by
$$\rho_2(g)[c] = [g_{11}][c]$$
and the representation 
$\sigma:P\by{} GL(V)$ is given by
$$\sigma(g)\left[\begin{array}{l}
b_1\\
\vdots\\ 
b_n\end{array}\right] =
[B]
\left[\begin{array}{l}
b_1\\
\vdots\\ 
b_n\end{array}\right]$$
which is the standard representation $\sigma:GL(n)\by{} GL(V)$.
Thus $$ \sT_{\bP_k^n}(-1) = \V_{\sigma} $$
is the  homogeneous bundle on $G/P$ associated to the standard
representation $\sigma$.
One can easily check that the maps $f$ and $\eta$ are compatible with
the $P$-module
structure of $V_2$, $V_1$ and $V$.

We write the sequence ~(\ref{e2}) as
$$0\by{} \V_{\rho_2}\by{} \V_{\rho_1}\by{} \V_{\sigma}\by{} 0.$$ 

\subsection{Restriction of $\V_{\sigma}$ to the quadric $Q\subset
\bP_k^n$}

  The bundle $\V_{\sigma} = \sT_{\bP_k^n}(-1)$, when restricted to $Q$,
fits into an  extension
\begin{equation}\label{e3}
0\by{} \sT_Q(-1) \by{} \sT_{\bP_k^n}(-1)\tensor_{\sO_{\bP_k^n}}\sO_Q
\by{} \sN_{Q/\bP_k^n}(-1)\by{} 0,
\end{equation}
where $\sT_Q$ and $\sN_{Q/\bP_k^n}$ denote the tangent sheaf and the 
normal sheaf
of
$Q\subset \bP_k^n$. Note that this is also a short exact sequence
of $SO(n+1)$-homogeneous bundles on
$Q = SO(n+1)/P_1$ (see \ref{n2}), hence there exists the corresponding  
short exact sequence of $P_1$-modules 
\begin{equation}\label{e4}
0\by{} U_1 \by{\tilde f} V \by{\tilde g}
U_{3}
\by{}
0, \end{equation}
where  $U_1$ and $U_3$ are $k$-vector spaces of dimensions
$n-1$ and $1$ respectively. We define
$${\tilde f}: (b_1,\ldots, b_{n-1})\to (b_1,\ldots, b_{n-1}, 0)$$
and 
$${\tilde g}: (a_1,\ldots, a_n)\to (a_n).$$ 
Now any $g\in P_1$ can be written as 
$$ g = \left[\begin{array}{lcl}
a_{11} &  * & *\\
0 & A & *\\
0 & 0 & a_{11}^{-1}\end{array}\right] $$

where $A\in
SO(n-1)$ and $a_{11}\in k\setminus \{0\}$.
 The representation
${\tilde\sigma}:P_1\by{} GL(V)$
is given by 
\begin{equation}\label{e8}
{\tilde\sigma}(g)\left[\begin{array}{l}
b_1\\ \vdots \\ b_n\end{array}\right]=
\left[\begin{array}{cc} A & *\\
                         0 &
a_{11}^{-1}\end{array}\right]\left[\begin{array}{l}
b_1\\ \vdots \\ b_n\end{array}\right]\end{equation}
The representation $\rho_3:P_1\by{} GL(U_3)$ is given by
$$\rho_3(g)[x] = [a_{11}^{-1}][x]$$
 and the representation $\sigma_1:P_1\by{} GL(U_1)$ is given by 

\begin{equation}\label{e9}
{\sigma_1}(g)\left[\begin{array}{l}
c_1\\ \vdots \\ c_{n-1}\end{array}\right] =
[A]\left[\begin{array}{l}
c_1\\ \vdots \\ c_{n-1}\end{array}\right]
\end{equation}
We write the sequence~(\ref{e4})
 as 
$$0\by{} \U_1\by{}\V_{\tilde \sigma} \by{} \U_3\by{} 0.$$

\begin{remark}\label{r1} Note that  $\sigma_1:P_1\by{} GL(U_{1})$ factors
through the standard
representation ${\tilde \sigma_1}:SO(n-1)\by{} GL(U_{1})$
and hence is irreducible, for $n\neq 3 $. This implies that the 
tangent bundle $T_Q$ is semistable. 
For $n=3$, the representation $\sigma_1$ is not irreducible and
$U_1$ is a direct sum of two $P_1$-submodules,
namely $k(1,0,0) \subset V$ and $k(0,1,0)\subset V$ respectively. In fact
one can check
easily that the only $P_1$-submodules of $V$ are given by
$k(1,0,0)$,
$k(0,1,0)$, $U_1$ and $V$ itself. In particular, 
all the homogeneous subbundles of $\V_{\tilde
\sigma}$ are given by these four $P_1$-submodules.

A smooth quadric $Q\subset
\bP^3_k$ is isomorphic to $\bP^1_k\times \bP^1_k$ and therefore the 
tangent bundle $T_Q$ is a direct sum of line
bundles of same degree. Hence the tangent bundle
$\sT_Q$ is always a semistable vector bundle for a smooth
 quadric $Q$.
Moreover, by
(\ref{e3}),
one can compute that $\mu(\sT_Q) >0$, if $n\geq 2$. 
\end{remark}

\section{Stablity  of  $\sT_{\bP_k^n}\mid_{\mbox{smooth quadric}}$}

\begin{proposition}\label{2}Let $\sigma:GL(n) \by{}
GL(V)$ be the standard representation
({\em i.e.}, $\sigma(g) = g$). Let $\V_{\sigma}$ be the
associated $G$-homogeneous bundle on $G/P = \bP_k^n$. Then
for characteristic $k \neq 2$, the restriction of the bundle $\V_{\sigma}$
= $\sT_{\bP_k^n}(-1)$
 to any smooth quadric $Q\subset \bP_k^n$ is
semistable.\end{proposition}

\noindent{\bf Remark} This result in characteristic $0$ is proved by
[F]. In
fact later we prove
a
stronger version of the above proposition (see Proposition~\ref{p3}). 

For the proof of the proposition we need the following two lemmas.

\begin{lemma}\label{l1}Let $\U_1$ and $\V_{\tilde \sigma}$
denote the $SO(n+1)$-homogeneous bundles, associated to the
$\sigma_1$ and
${\tilde \sigma}$ respectively (as given in Section~$2$), on $Q =
SO(n+1)/P_1$. Then 
$$\mu(\U_1) < \mu(\V_{\tilde \sigma}).$$\end{lemma}
\begin{proof}We are given that
$$\V_{\tilde \sigma} =
\V_{\sigma}\mid_Q~~ = \sT_{\bP_k^n}(-1)\tensor_{\sO_{\bP_k^n}}\sO_Q$$
and $\U_1 = \sT_Q(-1)$. Now

$$\deg~\sT_{\bP_k^n}(-1)\tensor_{\sO_{\bP_k^n}}\sO_Q
= 2~\deg~\sT_{\bP_k^n}(-1) = 2(\deg~H^0({\bP_k^n}, \sO_{\bP_k^n})
\tensor\sO_{\bP_k^n} -\deg~\sO_{\bP_k^n}(-1)) = 2,$$
where the second last equality follows from (\ref{e2}).
As 
$$\sN_{Q/\bP_k^n} \simeq (\sI/\sI^2)^{\vee} =
\sO_{\bP_k^n}(-2)^{\vee}\mid_Q =  \sO_{\bP_k^n}(2)\mid_Q, $$
where $\sI$ is the ideal sheaf of $Q\subset {\bP_k^n}$,
we have 
 $$\deg~\sN_{Q/\bP_k^n}(-1) = \deg~\sO_{\bP_k^n}(1)\mid_Q = 2.$$
Therefore 
$$\deg~\U_1 = \deg~\sT_Q(-1) = \deg~\sT_{\bP_k^n}(-1) -
\deg~\sN_{Q/\bP_k^n}(-1) =
0.$$
Hence $\mu(\U_1) = 0 < \mu(\V_{\tilde\sigma}) = 2/n $.
This proves the lemma.\end{proof}

\begin{lemma}\label{l2} The sequence~$(\ref{e4})$ 
$$0\by{} U_1 \by{\tilde f} V \by{\tilde g}
U_{3}
\by{}
0,$$
defined as above, 
of $P_1$-representations does not split.\end{lemma}
\begin{proof}It is enough to prove that the short exact
sequence~(\ref{e3}) does not split as sheaves of $\sO_Q$-modules. 
Suppose it does, then so does 
$$0\by{} \sT_Q(-2) \by{}
\sT_{\bP_k^n}(-2)\tensor_{\sO_{\bP_k^n}}\sO_Q
\by{} \sN_{Q/\bP_k^n}(-2) \by{} 0,$$
where we know that $\sN_{Q/\bP_k^n}(-2)\simeq \sO_Q$.
This implies that $H^0(Q, \sT_{\bP_k^n}(-2)\tensor_{\sO_{\bP_k^n}}\sO_Q) 
\neq 0$. However we have 
\begin{equation}\label{e11}
0\by{} \sT_{\bP_k^n}(-4) \by{} \sT_{\bP_k^n}(-2) 
\by{} \sT_{\bP_k^n}(-2)\tensor_{\sO_{\bP_k^n}}\sO_Q \by{} 0,
\end{equation}
where the first map is multiplication by the quadratic
equation defining
$Q\subset \bP_k^n$. If we assume the following   

\noindent{\bf Claim}.\quad $H^0(\bP_k^n, \sT_{\bP_k^n}(-2)) = 0 =
H^1(\bP_k^n,
\sT_{\bP_k^n}(-4))$.

Then (\ref{e11}) implies that 
$H^0(Q, \sT_{\bP_k^n}(-2)\tensor_{\sO_{\bP_k^n}}\sO_Q)= 0$, which
contradicts the hypothesis.
Now we give the 

\noindent{\underline{Proof of the claim}}. Consider the following 
 short exact sequence (which is derived from
(\ref{e2}))  
$$0\by{} \sO_{\bP_k^n}(-2) \by{} \sO_{\bP_k^n}(-1)^{n+1} \by{}
\sT_{\bP_k^n}(-2)\by{} 0. $$
As $n\geq 2$, we have $H^1({\bP}_k^n, \sO_{\bP_k^n}(-2))  = H^0({\bP}_k^n,
\sO_{\bP_k^n}(-1)) = 0 $, which implies 
$H^0(\bP_k^n, \sT_{\bP_k^n}(-2)) = 0$.
The above sequence also gives the long exact sequence 
$$ \by{}\oplus^{n+1}H^1(\bP_k^n,\sO_{\bP_k^n}(-3))\by{}
H^1(\bP_k^n,\sT_{\bP_k^n}(-4))\by{}
H^2(\bP_k^n,\sO_{\bP_k^n}(-4))\by{} $$
$$\by{}
\oplus^{n+1}H^2(\bP_k^n,\sO_{\bP_k^n}(-3))\by{} $$
\begin{enumerate}
\item If $n\geq 3$ then 
$H^1(\bP_k^n,\sO_{\bP_k^n}(-3))$ $ = H^2(\bP_k^n,\sO_{\bP_k^n}(-4))= 0$,
which implies $H^1(\bP_k^n,\sT_{\bP_k^n}(-4)) = 0$.
\item If $n=2$ then $H^1(\bP_k^2,\sO_{\bP_k^2}(-3)) = 0$. Moreover the map
$$H^2(\bP_k^2,\sO_{\bP_k^2}(-4))\by{}
\oplus^3H^2(\bP_k^2,\sO_{\bP_k^2}(-3)) $$
is dual to 
$$\oplus^3H^0(\bP_k^2,\sO_{\bP_k^2}) \by{}
H^0(\bP_k^2,\sO_{\bP_k^2}(1)) $$
which is an isomorphism as it comes from the evaluation map
$$H^0(\bP_k^2,\sO_{\bP_k^2}(1))\tensor\sO_{\bP_k^2} \by{}
\sO_{\bP_k^2}(1).$$
This implies
$H^1(\bP_k^n,\sT_{\bP_k^n}(-4)) = 0$.
\end{enumerate} 
This proves the claim and hence the lemma.\end{proof}

\vspace{10pt}

\noindent{\underline {\it Proof of Proposition}}~\ref{2}. Now suppose the
$SO(n+1)$-homogeneous bundle
$\V_{\tilde\sigma}$ on $Q$
 is not semistable.
Then it has a
Harder-Narasimhan filtration 
$$0\subset \V_1\subset \cdots \subset \V_{k} = \V_{\tilde\sigma} $$
where $\mu(\V_1) > \mu(\V_{\tilde\sigma})$. Now the uniqueness of the HN
filtration implies that $\V_1$ is a $SO(n+1)$-homogeneous subbundle of 
$\V_{\tilde\sigma}$. Therefore there exists a corresponding  
$P_1$-representation, say, 
$\rho_4:P_1\by{} GL({\tilde V}_1)$ and an inclusion of  $P_1$-modules
${\tilde V}_1 \hookrightarrow V$ corresponding
 to the
inclusion 
$\V_1\hookrightarrow \V_{\tilde\sigma}$.

\noindent{\bf Claim}. $U_1\subset {\tilde V}_1$, where $\sigma_1:P_1\to
GL(U_1)$ is
the $P_1$-representation as defined in (\ref{e9}).

We  assume the claim for the moment. Since 
$V/ U_1$ is an irreducible $P_1$-module, we have either
${\tilde V}_1 = U_1$ or 
$ {\tilde V}_1= V$, {\em i.e.}, $\V_1 = \U_1$ or
$\V_1 = \V_{\tilde\sigma}$. By Lemma~\ref{l1}, in both the cases
$\mu(\V_1) \leq
\mu(\V_{\tilde \sigma})$, which contradicts the fact that $\V_1$ is a term
of the HN filtration of $\V_{\tilde \sigma}$. Hence we conclude that the
$\V_{\tilde \sigma}$
is semistable.

Now we give 

\noindent{\underline{Proof of the claim}}. 
 Suppose ${\tilde V}_1\cap U_1 = 0$.
Then the composition map
$${\tilde V}_1 = \frac{{\tilde V}_1}{{\tilde V}_1\cap U_1}
\hookrightarrow \frac{V}{U_1} \hookrightarrow U_3,$$
gives an isomorphism ${\tilde V_1}\by{} U_3$, which implies that
(\ref{e4}) splits as
a sequence of $P_1$-modules; by Lemma~\ref{l2}, this is a
contradiction. 

Hence ${\tilde V}_1\cap U_1 \neq  0$.
If $n\neq 3$ then $U_1$ is  an
irreducible $P_1$-module (see Remark~\ref{r1}), which  implies that 
$U_1\subset
{\tilde V}_1$.  
Let $n = 3$ and $U_1\not\subset {\tilde V}_1$. Then 
Remark~\ref{r1} implies that $V_1\subset U_1$ as a $P_1$-submodule of rank
$1$ and 
therefore
$\mu(\V_1) =\mu(\U_1) <\mu(\V_{\tilde \sigma})$, which is a contradiction. 
 Therefore  $U_1\subseteq
{\tilde V}_1$. Hence the claim. This proves the proposition.
\hfill $\Box$

\begin{remark}\label{r2}The argument in the above proposition implies that
the only $SO(n+1)$-homogeneous subbundle of
$\sT_{\bP_k^n}(-1)\mid_Q$ $= \V_{\tilde \sigma}$ is 
either $\U_1$ or  
$\V_{\tilde \sigma}$ itself, if $n\neq 3$. If $n =3$ then the
homogeneous subbundle of $\V_{\tilde\sigma}$ is one of the two
homogeneous line subbundles of
$\U_1$ (as
given in Remark~\ref{r1}) or   
$\U_1$ or $\V_{\tilde \sigma}$ itself.\end{remark}

\begin{remark}\label{r5}For $n=3$, we can give another proof of  the
stability
of $\V_{\tilde \sigma}$ by reversing the role of cubic and quadric
in the proof of Lemma~\ref{l8}.   
\end{remark}

Now we can strengthen Proposition~\ref{2} as follows.

\begin{proposition}\label{p3}With the notations as in Proposition~\ref{2},
for $n\geq 3$, the restriction of the $\bP_k^n$-bundle,
$\V_{\sigma}$ to any smooth quadric $Q\subset \bP^n_k$ is
stable.
If $n=2$ then $\V_{\sigma}\mid_Q$ is a direct sum of two
copies of a line bundle on $Q$. \end{proposition}

Before coming to the proof of this proposition we need the following lemma
(which, perhaps, is already known to
the
experts). For this we recall some general facts. 
Let
$H$ be a reductive algebraic
group over $k$ and $P'\subset H$ be a parabolic group. 
Let $\V_{\rho}$ be a homogeneous $H$-bundle on $X=H/P'$ induced by 
a $P'$-representation $\rho: P'\by{} GL(V)$ on a
$k$-vector space $V$.  
Let the $H$ action on $\V_{\rho}$ be given by 
the map $L:H\times \V_{\rho}\by{} \V_{\rho}$, where
we write
$L(g, v) = L_g(v)$, for $g\in H$ and $v\in \V_{\rho}$.
This 
induces the canonical $H$-action on the dual of 
$\V_{\rho}$, which makes $\V_{\rho}^{\vee}$ and
$\V_{\rho}\tensor\V_{\rho}^{\vee}$
 into $H$-homogeneous
bundles
such that 
the map
$$\begin{array}{ccc}
\sEnd_{\sO_X}(\V_{\rho}) \tensor \sEnd_{\sO_X}(\V_{\rho}) &
\by{} &
\sEnd_{\sO_X}(\V_{\rho})\\
\downarrow_{\simeq} & & \downarrow_{\simeq}\\
(\V_{\rho}\tensor_{\sO_X} \V_{\rho}^{\vee})\tensor_{\sO_X}
(\V_{\rho}\tensor_{\sO_X} \V_{\rho}^{\vee}) & \by{} & 
(\V_{\rho}
\tensor_{\sO_X} \V_{\rho}^{\vee}).\end{array}$$
given by 
$$(v_1 \tensor \phi_1)\tensor(v_2\tensor \phi_2) \mapsto
\phi_1(v_2)(v_1\tensor\phi_2).$$ 
is $H$-equivariant. Hence 
$\End_{\sO_X}(\V_{\rho})  = 
H^0(X,\sEnd_{\sO_X}(\V_{\rho}))$ 
is a $H$-module such
that $H$
respects the algebra structure on it. This gives the homomorphism
$${\bar L}:H\by{} \mbox{Aut}(\End_{\sO_X}(\V_{\rho})),$$
given by 
${\bar L}(g)(\phi) = L_g\cdot\phi\cdot L_{g^{-1}},$
where $$\mbox{Aut}(\End_{\sO_X}(\V_{\rho})) = ~~\mbox{the set of ring
automorphism
on}~~\End_{\sO_X}(\V_{\rho}).$$ 

\begin{lemma}\label{l4}With the above notations, assume that the map
${\bar L}$,
defined as above, is
the trivial map. Then any subbundle of
$\V_{\rho}$ on $X$, which is also a direct summand of 
$\V_{\rho}$, is  $H$-homogeneous
vector subbundle.\end{lemma}
\begin{proof}Now let $\V_{\rho}= \U^1 \oplus \U^2$ be the direct sum of
subbundles $\U^1$ and $\U^2$. Let 
$\phi\in \End_{\sO_X}(\V_{\rho})$
 be given by 
  $$\phi\mid_{\U_1} = \rm{Id} ~~\mbox{and}~~\phi\mid_{\U_2} = 0.$$
Now, since ${\bar L}$ is 
trivial, we have
$${\bar L}(g)(\phi) = \phi ~~~\mbox{for all}~~~g\in G
.$$
{\it i.e.},
\begin{equation}\label{es}
L_g\cdot\phi\cdot L_{g^{-1}} = \phi. \end{equation} 
Let 
$(\V_{\rho})_x$ be the fiber of $\V_{\rho}$ over
$x\in
X$. 
Then, by (\ref{es}), we have the following commutative diagram
$$\begin{array}{ccc}
(\V_{\rho})_x & \by{L_{g^{-1}}} & (\V_{\rho})_{g^{-1}x}\\
\downarrow_{\phi} &  & \downarrow_{{\phi}_{g^{-1}}}\\
(\V_{\rho})_x & \by{L_{g^{-1}}} & (\V_{\rho})_{g^{-1}x},\end{array}$$
for each $x\in X$. This may be written as 
$$\begin{array}{ccc}
\U^1_x\oplus \U^2_x & \by{L_{g^{-1}}}
& \U^1_{g^{-1}x}\oplus\U^2_{g^{-1}x}\\
\downarrow_{\phi_x} &  & \downarrow_{\phi_{g^{-1}x}}\\
\U^1_x\oplus \U^2_x & \by{L_{g^{-1}}} & \U^1_{g^{-1}x}\oplus
\U^2_{g^{-1}x}.\end{array}$$
Now $$\U^2_x\subseteq \ker\phi_x \implies \U^2_x\subseteq 
\ker(L_g\cdot
\phi_{g^{-1}x}\cdot L_{g^{-1}}) = \ker(\phi_{g^{-1}x}\cdot L_{g^{-1}}) .$$ 
This implies 
$$L_{g^{-1}}(\U^2_x)\subseteq 
\ker~\phi_{g^{-1}x} = \U^2_{g^{-1}x}.$$
 Hence $L_{g^{-1}}(\U^2)\subseteq
\U^2$, {\it i.e.}, $\U^2$ is a $H$-homogeneous subbundle of $\V_{\rho}$.
This proves the lemma.\end{proof}

\vspace{5pt}

\noindent{\underline {\it Proof of Proposition}}~\ref{p3}.\quad By
Proposition~\ref{2},
for a quadric $Q\subset
\mathbf{P}_k^n$,
the bundle $\V_{\sigma}\mid_Q \simeq V_{\tilde\sigma}$ is
semistable. Hence there exists a nontrivial socle $\sF\subseteq
\V_{\tilde\sigma}$ such that $\mu(\sF) = \mu(\V_{\tilde\sigma})$ and $\sF$
is the maximal polystable subsheaf. Hence, by the uniqueness of maximal
polystable sheaf, it follows that it is 
an $SO(n+1)$-homogeneous 
subbundle of $\V_{\tilde\sigma}$. Therefore, by
Remark~\ref{r2}, either $\sF = \U_1$ or $\sF = \V_{\tilde\sigma}$. But 
$\mu(\sF) = \mu(\V_{\tilde\sigma}) > \mu(\U_1)$, which implies $\sF =
\V_{\tilde\sigma}$. Therefore  we can write 
$$\V_{\tilde\sigma} = \sF_1 \oplus \sF_2 \oplus \cdots \oplus \sF_r,$$
where $\sF_i$ is a direct sum of isomorphic stable sheaves, and the
stable summands of distinct $\sF_i$ are non-isomorphic. But
 each $\sF_i$ is an $SO(n+1)$-homogeneous subbundle of 
$\V_{\tilde\sigma}$ and is of the same slope as of
$\V_{\tilde\sigma}$. Hence
$r=1$ and $\V_{\tilde\sigma}$ is a direct sum of isomorphic stable
sub-bundles, {\it i.e.}
$$\V_{\tilde\sigma} = \oplus^{t}\U,
~~\mbox{where}~~\mu(\U) = \mu(\V_{\tilde \sigma}).$$ 
By Equation~(\ref{e2}), we have 
$$ 2 = \deg~\V_{\tilde\sigma} = t\cdot \deg~\U.$$
Hence $t =1 $ or $t =2$.

Suppose $n=2$. Then  $Q \simeq \bP^1_k$,
hence $\V_{\tilde\sigma}$ being rank~$2$ vector bundle on $Q$ splits as a
direct sum of two line bundles. Therefore in this case $t=2$. 

Suppose $n\geq 3$. If $t=1$ then we are done. 
Let $t =2$. Let 
$${\bar L}:SO(n+1)\by{}
\mbox{Aut}(H^0(Q,\sEnd(\V_{\tilde\sigma})))$$
 be the induced map. 
We are given that $\V_{\tilde\sigma} = \U\oplus\U$, where $\U$ is a stable
bundle on $Q$. But $\End_Q(\U)$ consists of scalars, and so 
$$\End_Q(\V_{\tilde\sigma}) \simeq M(2,k) ~~~\mbox{is
 the algebra of}~~~2\times 2
~~~\mbox{matrices}.$$
 Hence $\mbox{Aut}(H^0(Q,\sEnd(\V_{\tilde\sigma}))) \simeq
SO(3)$.
So, we have the map 
$${\bar L}: SO(n+1)\by{} SO(3).$$
But $SO(n+1)$ is an almost simple group, which implies, 
that 
$$\mbox{either}~~~\dim\mbox{Im}~{\bar L} = 0 ~~\mbox{or}~~~ 
\dim SO(n+1) = \dim\mbox{Im}~{\bar L} \leq \dim SO(3).$$
 Hence, for $n\geq
3$, $\dim\mbox{Im}~{\bar L}= 0$, which means ${\bar L}$ is trivial. 
Therefore, by
Lemma~\ref{l4}, the bundle $\U$ is homogeneous.

However, by Remark~\ref{r2} and Lemma~\ref{l1}, the only $G$-homogeneous
subbundle of
$\V_{\tilde \sigma}$, of the same slope as  $\V_{\tilde \sigma}$, is
$\V_{\tilde \sigma}$ itself. 
Hence we conclude that $\V_{\tilde \sigma} = \U$ is stable, if $n\geq 3$. 
This proves the proposition.\hfill $\Box$ 

\begin{cor}\label{c3}If $Q\subset \bP^n_k$ is a smooth quadric such that
$k$ is an algebraically closed field of char $\neq 2$ then 
\begin{enumerate}
\item $\Omega_{\bP^n_k}\mid_Q$ is strongly semistable if $n=2$ and
\item $\Omega_{\bP^n_k}\mid_Q$ is strongly stable if $n\geq
3$.\end{enumerate}
\end{cor}
\begin{proof}If $n =2$ then the corollory follows from
Proposition~\ref{p3}. Suppose  $n\geq 3$. Then, by Proposition~\ref{p3},
the bundle $\Omega_{\bP^n_k}\mid_Q$ is stable. Moreover, by
Remark~\ref{r1}, the tangent bundle $\sT_Q$ of $Q$  is semistable and 
$\mu(\sT_Q)>0$. Hence, by Theorem~2.1 of [MR1], the bundle 
$\Omega_{\bP^n_k}\mid_Q$
is strongly stable. This proves the corollory.\end{proof}

\section{Stablity  of  $\sT_{\bP_k^n}\mid_{\mbox{smooth cubic}}$}   
We recall the Bott vanishing theorem for $(\bP^n_k,
\Omega^q_{\bP^n_k}(t))$, where $k$ an
arbitrary field of arbitrary characteristic.
$$\begin{array}{lcl}
H^0(\bP^n_k, \Omega^q_{\bP^n_k}(t)) & \neq & 0, ~~\mbox{if}~~0\leq q\leq
n,~~~\mbox{and}~~~t > q\\   
H^n(\bP^n_k, \Omega^q_{\bP^n_k}(t)) & \neq  & 0 ~~\mbox{if}~~~0\leq q\leq
n,~~~\mbox{and}~~t < q-n\\
H^p(\bP^n_k, \Omega^p_{\bP^n_k}) & = & k,~~~\mbox{if}~~0\leq p\leq n\\
H^p(\bP^n_k, \Omega^q_{\bP^n_k}(t)) & = & 0~~~\mbox{otherwise}.
\end{array}$$

Now throughout this section we fix a smooth hypersurface $X$ of degree
$d\geq 3$ in
$Y = \bP^n$, $(d, \mbox{char}~k) =1$.
We have the following short exact sequences
\begin{equation}\label{e17}
0\by{} \Omega^q_Y(t) \by{} \Omega^q_Y(t+d) \by{}
\Omega^q_Y(t+d)\mid_X\by{} 0\end{equation}

\begin{equation}\label{e18}
0\by{} \Omega^q_X(t) \by{} \Omega^{q+1}_Y(t+d)\mid_X \by{}
\Omega^{q+1}_X(t+d)\by{} 0\end{equation}

\begin{enumerate}
\item If $p+q < \dim~X$ and $p, q\geq 0$ then from  Bott vanishing
 and the  short exact sequences  (\ref{e17}) and  (\ref{e18}), it follows
that  
$H^p(X,\Omega^q_X(t))= 0$ for $t<0$. 
\item If $p+q < \dim~X$ then
$$H^p(X, \Omega^q_X) \simeq H^p(Y,\Omega^q_Y).$$
\item Consider the following commutative diagram of natural maps
$$\begin{array}{ccc}
H^{p}(Y,\Omega^{q}_Y) & \by{} & H^{p+1}(Y,\Omega^{q+1}_Y)\\
\downarrow & & \downarrow\\
H^{p}(X,\Omega^{q}_X) & \by{} & H^{p+1}(X,\Omega^{q+1}_X),
\end{array}$$
where the horizontal maps are given by the cup product with
$c_1(\sO_Y(d))   
= d\cdot c_1(\sO_Y(1))$ and $c_1(\sO_X(d))$ respectively. 
Since $(\rm{char}~k, d) =1$, the
map $H^p(Y, \Omega^q_Y)\by{} H^{p+1}(Y, \Omega^{q+1}_Y)$ is an isomorphism
for
every $p,q$ with $p,q\geq 0$ and $p+1\leq \dim~Y$. 
In particular, the induced composite map
\begin{equation}\label{e19}
\eta_{p,q}:H^{p}(X,\Omega^{q}_X)\by{} H^{p+1}(Y,\Omega^{q+1}_Y)
\end{equation}
is an isomorphism if $p,q\geq 0$ and $p+q < \dim~X$.
\end{enumerate}

We prove the following Lemma~\ref{l5} and Corollory~\ref{c2} along the
same line of arguements, as 
given  for the case $k =\C$, in
[PW].
\begin{lemma}\label{l5}Let $X\subseteq \bP^n_k$ be a hypersurface of
$\deg~d\geq 3$. Let $n\geq 2$ and $(\rm{char}~k, d) = 1$. If $p, q
\geq 0$ and $p+q <
\dim~X$
and $t\leq q(n+1-d)/(n-1)$ then
\begin{enumerate}
\item $H^p(X, \Omega^q_X(t)) = 0$, if $t\neq 0$ and
\item $H^p(X, \Omega^q_X) \simeq H^p(Y,\Omega^q_Y) $.\end{enumerate}
\end{lemma}
\begin{proof}As discussed above, (a)~for $t<0$, the statement~(1)  holds,
 {\it i.e.}, for $t<0$, we have $H^p(X,
\Omega^q_X(t)) =0$,  and (b)~the statement~(2) always holds.
 
 Suppose $t =d$. In particular $q\geq 2$.
 Now
$(\ref{e18})$ gives the long
exact sequence
$$H^p(\Omega^{q-1}_X)\by{f_{p,q-1}} H^p(\Omega^q_Y(d)\mid_X)\by{}
H^p(\Omega^{q}_X(d))\by{}
H^{p+1}(\Omega^{q-1}_X)\by{f_{p+1,q-1}} H^{p+1}(\Omega^{q}_Y(d)\mid_X).$$
Hence to prove that $H^p(X, \Omega^q_X(d)) = 0$, it is enough to prove the
following

\noindent{\underline{Claim}}: The map $f_{p,q}$ is an
isomorphism, if $p, q\geq 0$ and 
$p+q<\dim~X$.

\vspace{5pt}  

\noindent{\underline{Proof of the claim}}. Note that we have the following
commutative diagram
$$\begin{array}{lcl}
H^{p}(X,\Omega^q_X) & \by{f_{p,q}} & H^{p}(Y, \Omega^{q+1}_Y(d)\mid_X)\\
 {} & \searrow^{\eta_{p,q}} & \downarrow^{g_{p,q+1}}\\
{} & {} & H^{p+1}(Y,\Omega^{q+1}_Y),\end{array}$$
where, by $(\ref{e19})$ the map $\eta_{p,q}$ is an isomorphism.
 Hence the map $g_{p,q+1}$ is surjective, in this
case. Moreover, by $(\ref{e17})$  we also have the exact sequence

$$H^p(Y, \Omega_Y^{q+1}(d)) \by{}
H^p(X, \Omega_Y^{q+1}(d)\mid_X) \by{g_{p,q+1}}
H^{p+1}(Y, \Omega^{q+1}_Y),$$

where $H^p(Y, \Omega_Y^{q+1}(d)) = 0$, by Bott vanishing. Therefore the
map
$g_{p,q+1}$ is an isomorphism. This implies that
$f_{p,q}$ is an isomorphism.
 This proves the claim.   
Hence $H^p(X, \Omega^q_X(d)) =0$ if $p,q \geq 0$ and $p+q < \dim~X$
\vspace{5pt}

By induction on $t$, we can assume that for $m < t$ and $m\neq 0$, 
 we have 
$$H^i(X, \Omega^j_X(m)) = 0, ~~\mbox{where}~~ i,j\geq 0, ~~i+j < 
\dim~X 
~~~\mbox{and}~~
 m\leq \frac{j(n+1-d)}{n-1},$$ 

Now, to prove the proposition, it remains to show that, 
$$t \leq \frac{q(n+1-d)}{(n-1)},~~t\not\in \{0,d\}, ~~p, q \geq
0, ~~p+q <\dim~X \implies
H^p(X, \Omega_X^q(t)) = 0.$$

Note that the hypothesis that 
$$t \leq \frac{q(n+1-d)}{n-1} \implies t\leq q.$$ 

Consider the following long exact sequence (obtained from $(\ref{e18})$)

$$ H^p(X, \Omega_Y^q(t)\mid_X)\by{}
H^p(X, \Omega_X^{q}(t))\by{} H^{p+1}(X, \Omega^{q-1}_X(t-d))$$

If $q-1 < 0$ then the last term is $0$. If $q-1\geq 0$ then as 
$$t \leq \frac{q(n+1-d)}{n-1} \implies  t-d \leq 
\frac{(q-1)(n+1-d)}{n-1},$$
 by induction hypothesis on $t$, the last term
of the
sequence is $0$. Consider the exact sequence (obtained from
$(\ref{e17})$)
$$H^p(Y, \Omega_Y^{q}(t)) \by{} H^p(X, \Omega_Y^{q}(t)\mid_X)\by{}
H^{p+1}(Y, \Omega^{q}_Y(t-d))$$
then, by  Bott
vanishing, the first and the last term of the sequence are $0$. This
implies that   
$H^p(X, \Omega_Y^q(t)\mid_X) = 0$. Hence $H^p(X, \Omega_X^q(t)) =0$. This
completes the proof of the proposition.\end{proof}

\begin{cor}\label{c2}Let  $X\subset \bP^n_k$  be a smooth 
hypersurface
of
degree $d\geq 3$. Let $n\geq 4$ and $g.c.d.(\rm{char}~k, d) = 1$.
Then  $\Omega_X$ is stable.
\end{cor}
\begin{proof}Suppose $\Omega_X$ is not stable then there
exists a subbundle  $W\subset \Omega_X$ of rank $q\leq
n-2$, such that 
$\mu(W) \geq \mu(\Omega_X)$. 
Then $\wedge^qW\hookrightarrow \wedge^q\Omega_X$. Since $\wedge^qW\in
\Pic(X)$, we have $\wedge^qW = \sO_{\bP^n_K}(-t)\mid_X$,
 as $n\geq 4$ implies that the map  $\Pic(\bP^n_K)\to \Pic(X)$ is an
isomorphism. This implies that 
$H^0(X,\Omega_X(t)) \neq 0$.
Hence to prove that the bundle $\Omega_X$ is stable, it is
enough to prove that 
$$H^0(X,\Omega^q_X) = 0, ~~~\mbox{for}~~t\leq \frac{q(n+1-d)}{n-1},$$
which  immediately follows by Lemma~\ref{l5}. Hence $\Omega_X$ is stable. 
\end{proof}

\begin{lemma}\label{l7}Let $X\subset \bP^3_k$ be a smooth hypersurface
of degree $d =3$. Then $\mu_{\rm min}(\sT_X) \geq 0$.\end{lemma}
\begin{proof}Let $H\subset \bP^3_k$ be a general hyperplane such that $C =
X\cap
H$ is a nonsingular complete intersection on $\bP^3_k$. In particular $C$
is an elliptic curve. This gives the canonical  short exact sequence 
$$0\by{} \sT_C \by{} \sT_X\mid_C \by{} \sN_{C/X}\by{} 0, $$
which is equivalent to 
$$0\by{} \sO_C \by{f_1} \sT_X\mid_C \by{f_2} \sO_C(1)\by{} 0.$$
If $T_X$ is semistable then $\mu_{\rm min}(\sT_X) = \mu(\sT_X) = 1/2 >0$.
We can assume that  $T_X$ is not semistable.
Let $\sL \subset \sT_X$  be the Harder-Narasimhan filtration of $\sT_X$,
which gives a short exact sequence of coherent sheaves (where $\sL$ is a
line
bundle on $X$),
$$0\by{} \sL \by{g_1}\sT_X \by{g_2}\sM\by{} 0.$$
By definition,  $\mu_{\rm min}(\sT_X) = \deg~\sM$, therefore it is enough
to
prove that
$\deg~\sM >0$, which is same as to prove that 
$\deg~\sM\mid_C = \sM\cdot H >0$.
Consider the composite map 
$$\sO_C\by{f_1}\sT_X\mid_C\by{g_2\mid_C} \sM\mid_C.$$
{\underline {Case}}~{1}. If $g_2\mid_C \circ f_1 = 0$ then the induced map
$\sO_C(1)\by{}
\sM\mid_C$ is surjective. This implies that $\deg~\sM\mid_C>0$.
{\underline {Case}}~{2}. If $g_2\mid_C \circ f_1 \neq 0$ then there exists
a
nonzero map 
$\sO_C\by{} \sM\mid_C$, which implies that $\deg~\sM\mid_C\geq 0$. 
 This proves the
lemma.\end{proof} 

\begin{lemma}\label{l6}Let $X\subset \bP^n_k$  be
a smooth hypersurface of degree $d\geq 3$. Let $n\geq 4$ and
$g.c.d.({\rm char}~k, d) = 1$. 
Then $\Omega_{\bP^n_k}\mid_X$ is stable.
\end{lemma}
\begin{proof}As  argued in Corollory~\ref{c2}, it is
enough to 
 prove that 
$$H^0(X,
\Omega^q_{\bP^n_k}(t)\mid_X) = 0, ~~~\mbox{for}~~~t\leq
q(n+1)/n~~~\mbox{and}~~~1\leq q\leq n-1.$$
Now, consider 
$$0\by{} \sO_{\bP^n_k}(-d)\by{} \sO_{\bP^n_k}\by{} \sO_X\by{} 0,$$
which gives 
$$0\by{} \Omega^q_{\bP^n_k}(t-d)\by{} \Omega^q_{\bP^n_k}(t)\by{} 
\Omega^q_{\bP^n_k}(t)\mid_X\by{} 0.$$
Since $t\leq q(n+1)/n \implies t\leq q$, by Bott vanishing we have 
$$H^0(\bP^n_k, \Omega^q_{\bP^n_k}(t)) = 0, ~~~\mbox{for}~~ t\leq
q(n+1)/n,$$
and
$$H^1(\bP^n_k, \Omega^q_{\bP^n_k}(t-d))  = 0, ~~~\mbox{if}~~t\neq
d~~~\mbox{or}~~q\neq 1.$$
Therefore the exact sequence 
$$H^0(\bP^n_k, \Omega^q_{\bP^n_k}(t))\by{} 
H^0(\bP^n_k, \Omega^q_{\bP^n_k}(t)\mid_X)\by{}
H^1(\bP^n_k, \Omega^q_{\bP^n_k}(t-d))$$
implies that for $t\leq q(n+1)/n$
$$H^0(\bP^n_k, \Omega^q_{\bP^n_k}(t)\mid_X) = 0, ~~\mbox{if}~~t\neq
d~~~\mbox{or}~~q\neq 1.$$
However the case, when $t= d$ and $q=1$ and $t\leq q(n+1)/n$ does not
arise, as these conditions 
 imply that $d = t\leq 1+ (1/n) < 2$. Hence we conclude that 
$H^0(\bP^n_k, \Omega^q_{\bP^n_k}(t)\mid_X) = 0$ if $t\leq q(n+1)/n$. This
proves the lemma.\end{proof}

\begin{lemma}\label{l8}Let $X\subset \bP^n_k$ be a smooth cubic
hypersurface such that $n=2$ or $n=3$. Then $\Omega_{\bP^n_k}\mid_X$ is
strongly semistable.
\end{lemma}
\begin{proof}Suppose $n=2$, then $X$ is an elliptic curve. Hence 
$\Omega_{\bP^2_k}\mid_X$ is an indecomposable  rank $2$ vector bundle on
$X$ (see the proof of Theorem~3.16 of [NT]) and is of negative
degree.  Hence
strong semistabilty follows from the facts that a vector bundle of
negative
degree  has no sections and a semistable bundle is strongly
semistable on an elliptic curve.

Suppose $n=3$. Let  $Q\subset \bP^3_k$ be a general smooth quadric such
that $C =
Q\cap X$ is a smooth complete intersection nonsingular curve in
$\bP^3_k$. Then $C$ is  curve of genus $=4$ such that 
$\sO_{\bP^3_k}(1)\mid_C = \omega_C$ and the restriction of the short exact
sequence 
$$0\by{} \Omega_{\bP^3_k}(1)\by{} H^0(\bP^3_k, \sO_{\bP^3_k}(1))\tensor
\sO_{\bP^3_k} \by{} \sO_{\bP^3_k}(1)\by{} 0,$$
to  $C$, is 
$$0\by{} \Omega_{\bP^3_k}(1)\mid_C \by{} 
H^0(C,\omega_C)\tensor\sO_C \by{}
\omega_C\by{} 0.$$
Note that $C$ is a non-hyperelliptic curve, hence by Corollory~3.5 of
[PR] (the proof given there for $k = \C$ works for 
any algebraically closed field $k$ of arbitrary
characteristic), the bundle $\Omega_{\bP^3_k}(1)\mid_C$ is stable. By
Lemma~\ref{l7}, we have $\mu_{\rm min}(\sT_X) \geq 0$. Therefore
Theorem~2.1 of [MR1] implies that $\Omega_{\bP^3_k}(1)\mid_C$ is
strongly semistable, for general curve $C\subset X$, of degree $3$. Hence
$\Omega_{\bP^3_k}(1)\mid_X$ is strongly
semistable. Hence the lemma.
\end{proof}

\begin{cor}\label{c4}If $X\subset \bP^n_k$ is a smooth cubic such that $k$
is an algebraically closed field of characteristic~$\neq 3$, then 
\begin{enumerate}
\item $\Omega_{\bP^n_k}\mid_X$ is strongly semistable, if $n=2$
or $3$ and 
\item $\Omega_{\bP^n_k}\mid_X$ is strongly stable, if $n\geq 4$
\end{enumerate}
\end{cor}
\begin{proof}The cases $n=2$ and $n=3$ follow from
Lemma~\ref{l8}. Hence it is
enough to prove the corollory for $n\geq 4$. Now, by Corollory~\ref{c2},
the
tangent bundle $\sT_X = \Omega^{\vee}_X$ of $X$ is semistable and is of
positive slope. By
Lemma~\ref{l6}, the bundle $\Omega_{\bP^n_k}\mid_X$ is stable. Hence,
 again, by Theorem~2.1 of  [MR1], we deduce that
$\Omega_{\bP^n_k}\mid_X$ is strongly stable. Hence the
corollory.\end{proof}

 \section{Main results}

\begin{notation}\label{n3}We recall  the notion of `generic' and
`general' as
given in Section~1
of [MR2]. Let $k$ be an algebraically closed field of arbitrary
characteristic.
Let $S_d = {\rm Proj}(H^0(\bP^n_k, \sO_{\bP^n_k}))$. Then we have 
$$\begin{array}{ccc}
\bP^n_k\times S_d\supseteq Z_d & \by{q_d} & S_d\\
\downarrow^{p_d} & &\\
\bP^n_k & & ,\end{array}$$
wherw $Z_d = \{(x,s)\in \bP^n_k\times S_d\mid s(x) = 0\}$ and $p_d$, $q_d$
are projections. The fiber of $q_d$ over $s\in S_d$ is the embedding in
$\bP^n_k$ via $p_d$ as the hypersurface of $\bP^n_k$ defined by the ideal
generated by $s$. 
Let $K_d$ be the function field of $S_d$. Let $Y_d$ be the generic fiber
of $q_d$ given by the fiber product 
$$\begin{array}{ccc}
Z_d & \by{} & S_d\\
\uparrow^{q_d} & & \uparrow\\
Y_d & \by{} & {\rm Spec}~K_d,\end{array}$$
where $Y_d$ is an absolutely irreducible, nonsingular hypersurface, and
there
is a nonempty open subset of $S_d$ over which the geometric fibres of
$q_d$
are irreducible.

We call $Y_d$ {\it the generic hypersurface} of degree $d$. Whenever a
property holds for $q_d^{-1}(s)$ for $s$ in a nonempty Zariski open subset
of
$S_d$, then we say it holds for {\it a general} $s$. 
 \end{notation}

\begin{remark}\label{r4}For a torsion free sheaf $V$ on a smooth
projective variety (which is $\bP^n_k$ in our case), the restriction of
$V$
to the generic hypersurface $Y_d$ is semistable (geometrically stable) if
and only if the 
restriction of $V$ to a general hypersurface of degree $d$ is semistable
(geometrically stable): because, for any coherent torsion free sheaf $F$
of  $X$, the sheaf $p_d^*F$ forms a flat family over a nonempty open
subset
of $S_d$ (see Proposition~1.5 of [MR2]), and the property of coherent
sheaves being
semistable (geometrically stable) is open in flat families.\end{remark}

\begin{remark}\label{r6}If
\begin{enumerate} 
\item $X$ = smooth quadric, if ${\rm char}~k \neq 2$, or
\item $X$ = smooth  cubic, if ${\rm char}~k \neq 3$ 
\end{enumerate}
then, by Corollory~\ref{c3} and Corollory~\ref{c4}, the bundle
$\Omega_{\bP^n_k}\mid_X $ 
is
strongly semistable. Moreover, by Remark~\ref{r1}, Corollory~\ref{c2} and
Lemma~\ref{l7}, we have 
$\mu_{\rm min}(\sT_X) \geq 0$.
 In particular, by Theorem~2.1 of [MR1] and Theorem~3.23 of [RR], any
semistable bundle on $X$ remains semistable 
 after applying the  functors like Frobenius pull backs, tensor powers, 
symmetric powers,
and exterior powers on $X$.\end{remark}

\noindent{\underline {\it Proof of Theorem}}~\ref{1}. By Remark~\ref{r6},
it is
enough to prove that $\W_{\tau}$ is semistable on $X$. By
Proposition~2.4
of [J], given an irreducible representation 
$$\tau:GL(n)\by{} GL(W),$$ 
 there exists  $\lambda\in \chi(T)$ (for a
fixed torus $T$ of
$GL(n)$)
such that 
$$W = L(\lambda),$$
 where 
following the notation of [J], the $GL(n)$-module $L(\lambda)$ = socle of
$H^0(\lambda)$. Moreover, by corollory~2.5 of [J], the module dual to 
$L(\lambda)$ is 
$$L(\lambda)^{\vee} = L(-w_0\lambda).$$
Let $\epsilon_i \in \chi(T)$ be given by
$\epsilon_i(t_1,t_2,\ldots,t_n) = t_i$ and let $\omega_i =
\epsilon_1+\cdots +\epsilon_i$. Then      
 any $\nu\in \chi(T)$ can be
written as 
$$\nu = \sum_ia_i\omega_i = \sum_i\nu_i\epsilon_i,$$
where $\nu_i\in \Z$ and $\nu_1\geq \nu_2\geq \cdots \geq \nu_n$.

Let $\Hm^0(L_{\nu})$ be the vector bundle on $G/P =
\bP_k^n$ 
corresponding to the
$GL(n)$-representation $H^0(L_{\nu})$.

\noindent{\bf Claim}. The bundle $\Hm^0(L_{\nu})\mid_X$ is semistable on
$X\subset {\bP}_k^n$ and 
$$\mu(\Hm^0(L_{\nu})\mid_X) =
(\sum_i\nu_i)(\mu(\V_{\sigma}\mid_X)),$$
 
\noindent{\it Proof of the claim}: Let us  denote 
$$ S(a_1,\ldots,a_n, V) = S^{a_1}(V)\tensor
S^{a_2}(\wedge^2V)\tensor\cdots
\tensor S^{a_n}(\wedge^nV),$$
for a vector space $V$, and let us denote

$$S(a_1,\ldots,a_n, \V) = S^{a_1}(\V)\tensor
S^{a_2}(\wedge^2\V)\tensor\cdots
\tensor
S^{a_n}(\wedge^n\V),$$
for a vector bundle $\V$.
By definition of
$H^0(L_{\nu})$, we
have a surjection of $GL(n)$-modules 

\begin{equation}\label{e5}  S(a_1,\ldots,a_n, V)\by{} H^0(L_{\nu}),
\end{equation}
where $\sigma:GL(n)\by{} GL(n)=GL(V)$ is the standard representation.
Hence we have the surjection of $G$-homogeneous bundles on ${\bP}_k^n$ 
\begin{equation}\label{e10}
S(a_1,\ldots,a_n, \V_{\sigma}) \by{} \Hm^0(L_{\nu}),   
\end{equation}
where we recall that $\V_{\sigma}= \sT_{\bP^n_k}(-1) =
(\Omega_{\bP^n_k}(1))^{\vee}$ is the vector bundle associated to the
representation $\sigma$.
Therefore we have the surjection of bundles on $X$ 
\begin{equation}\label{e14}
S(a_1,\ldots,a_n, \V_{\sigma}\mid_X) \by{} \Hm^0(L_{\nu})\mid_X.   
\end{equation}

By Theorem~1.1 (and Cor.~1.3), expos$\acute{e}$ XXV, Sch$\acute{e}$mas en
groupes III,  
[SGA-3], $GL(n)/B$ ($B$ is a Borel group of $GL(n)$) can be lifted to
characteristic zero. Therefore
 the degree and rank of these vector bundles are independent of
the characteristic of the field. Now over a field of characteristic $0$,
sequence (\ref{e5}) split, which implies that sequence~(\ref{e10}) splits
as bundles on $\bP_k^n$, defined over field of characteristic $0$.
Now since $S(a_1, \ldots, {a_n}, \V_{\sigma})$
is semistable vector bundle, we have
$$\begin{array}{ccc}
\mu(\Hm^0(L_{\nu})) & = & \mu(S(a_1, \ldots, {a_n}, \V_{\sigma}))\\
 & = & (a_1 + 2a_2 + \cdots +na_n)\mu(\V_{\sigma})\\ 
 & = & (\sum_i\nu_i)\mu(\V_{\sigma}),\end{array}$$
where the last inequality follows as $\nu_i = a_i+\cdots a_n $.
Hence  

\begin{equation}\label{e12}\mu(\Hm^0(L_{\nu})\mid_X) =
(\sum_i\nu_i)(\mu(\V_{\sigma}\mid_X)).\end{equation}

 By Remark~\ref{r6}, the bundle $S(a_1, \ldots,
{a_n},
\V_{\sigma}\mid_X)$
is semistable. Therefore, by (\ref{e14}) and (\ref{e12}), 
the bundle $\Hm^0(L_{\nu})\mid_X$ is
semistable. Hence the claim.

Now, coming back to  $W = L(\lambda)$, let
$$\lambda = \sum_ia_i\omega_i = \sum_i\lambda_i\epsilon_i.$$
Then, as $w_0(\epsilon_i) = \epsilon_{n+1-i}$, we have
$$-w_0\lambda =
a_{n-1}\omega_1+\cdots + a_1\omega_{n-1} +(-a_1+ \cdots
-a_n)\omega_n = -\sum_i(\lambda_{n+1-i})\epsilon_i.$$
This implies that $\mu(\Hm^0(L_{-w_0\lambda})) =  
-\mu(\Hm^0(L_{\lambda}))$, therefore 

\begin{equation}\label{e16}\mu(\Hm^0(L_{-w_0\lambda})\mid_X) =  
-\mu(\Hm^0(L_{\lambda})\mid_X).\end{equation}
 
Moreover there exists the surjective map of vector
bundles on $X$

\begin{equation}\label{e7}
S(a_1, \ldots, {a_n}, \V_{\sigma}\mid_X)
\tensor S(a_{n-1}, \ldots, {a_1},-(a_1+\cdots+a_n), \V_{\sigma}\mid_X)
\by{} (\Hm^0(L_{\lambda})\tensor
\Hm^0(L_{-w_0\lambda}))\mid_X,
\end{equation}

where the L.H.S. is a semistable vector bundle  of slope
$= 0$. Moreover, by (\ref{e16}), the slope of
R.H.S. is also $ = 0$. Hence
$\Hm^0(L_{\lambda})\mid_Q\tensor
\Hm^0(L_{-w_0\lambda})\mid_X$ is semistable of slope $0$.
Now, consider the injective map
$$\W_{\tau}\tensor \W_{\tau}^{\vee} \by{}
\Hm^0(L_{\lambda})\tensor
\Hm^0(L_{-w_0\lambda}),$$
which give the injective map
\begin{equation}\label{e6}
\W_{\tau}\mid_X\tensor \W_{\tau}^{\vee}\mid_X \by{}
\Hm^0(L_{\lambda})\mid_X\tensor
\Hm^0(L_{-w_0\lambda})\mid_X
\end{equation}
is injective, where the slope of L.H.S is $=0$, which is same as the slope
of   R.H.S..
Hence
$\W_{\tau}\mid_X\tensor
\W_{\tau}^{\vee}\mid_X$ is semistable. This implies that $\W_{\tau}\mid_X$
is semistable, which proves the theorem.\hfill $\Box$

\begin{cor}\label{c1}Let $\W_{\tau}$ be the homogeneous bundle on
$\bP_k^n$ associated to an irreducible representation $\tau:GL(n)\by{}
GL(W)$. Let $k$ be an algebraically closed field of characteristic
$\neq 2, 3$. Then 
\begin{enumerate}
\item for $s\geq 0$, the $s^{th}$ Frobenius power
$F^{s*}\W_{\tau}\mid_H$ is
semistable, for general hypersurface $H$ of  degree $d\geq 2$ in
$\bP_k^n$. In particular
\item  $\W_{\tau}\mid_{H_0}$
is strongly semistable, where $H_0\subset \bP^n_{K_d}$ is the 
 $k$-generic
 hypersurface  of 
degree $d \geq  2$.\end{enumerate}
Moreover,
if $\W_{\tau}$ is the  tangent bundle on $\bP^n_k$ and $n\geq 4$ then
we
can replace the
word `semistable' by
`stable' everywhere in the above statement.
\end{cor}
\begin{proof}By Theorem~\ref{1}, the bundle
$\W_{\tau}\mid_X$ is strongly semistable, where
$X$ is a smooth quadric or a smooth cubic in $\bP^n_k$.
In other words, for $s\geq 0$ and for the $s^{th}$ iterated
Frobenius pull back, $F^{s*}\W_{\tau}$ of $\W_{\tau}$, 
the bundle $F^{s*}\W_{\tau}\mid_X$  is semistable, where
 $X$ is a smooth quadric or a smooth  cubic. 
 Hence, by the proof of the restriction
theorem of [MR2], it follows that $F^{s*}\W_{\tau}\mid_H$ is
semistable when restricted to a general hypersurface $H\subset \bP^n_k$ of
 degree $\geq 2$ (see also the modified proof of the above mentioned
restriction
theorem given in [HL]).
This proves  part~(1) of the corollory.

Moreover this implies that, for any $s\geq 0$ and for generic 
hypersurface $H_0$ of degree $\geq 2$, the bundle
$F^{s*}\W_{\tau}\mid_{H_0}$ is
semistable (see Remark~\ref{r4}). In particular, the bundle
$\W_{\tau}\mid_{H_0}$ is strongly
semistable. This
proves the
part~(2) of the corollory.

Note that, for $n\geq 4$, by Corollories~\ref{c3} and \ref{c4}, the 
bundle $\sT_{\bP_k^n}\mid_X$ is strongly  stable and hence
geometrically strongly stable (as the
underlying field $k$ is algebraically closed). Now the similar arguments,
as above, 
applied to the tangent
bundle $\sT_{\bP^n_k}$,  prove the rest of  the corollory.
\end{proof}

\begin{remark}\label{r3}By Proposition~\ref{p3}, the bundle
$\sT_{\bP_k^n}\mid_Q$ is stable for a smooth quadric $Q\subset
\bP_k^n$, for $n\geq 3$. One may ask the following: If $\tau:GL(n)\by{}
GL(W)$
is an irreducible representation, then is the associated 
bundle $\W_{\tau}$  stable  on $Q$? More generally if $\tau:GL(n)\by{}
H$ is any irreducible representation, with $H$ semisimple, then is the
induced $H$ bundle semistable on $Q$? \end{remark}

\end{document}